\documentclass[a4paper,12pt,reqno]{article}

\usepackage{amsmath,amsfonts,amsthm,amssymb}
\usepackage{isolatin1}

\newcommand{\dom}{\mbox{Dom}}

\newcommand{\iot}{\int_{0}^{t}}
\newcommand{\ist}{\int_{s}^{t}}
\newcommand{\iotd}{\int_{[0,t]^2}}
\newcommand{\iott}{\int_{0}^{T}}
\newcommand{\iou}{\int_{0}^{1}}
\newcommand{\kh}{K_H}

\newcommand{\me}{M^{\ep}}
\newcommand{\mes}{M^{\ep,*}}
\newcommand{\ms}{M^*}
\newcommand{\ot}{[0,t]}
\newcommand{\ott}{[0,T]}
\newcommand{\ou}{[0,1]}
\newcommand{\xe}{X^{\ep}}
\newcommand{\1}{{\bf 1}}

\newcommand{\R}{\mathbb R}

\newcommand{\cf}{\mathcal F}
\newcommand{\ch}{\mathcal H}

\newcommand{\al}{\alpha}
\newcommand{\ep}{\varepsilon}

\newcommand{\ga}{\gamma}
\newcommand{\ka}{\kappa}

\newcommand{\si}{\sigma}

\newcommand{\vp}{\varphi}

\newcommand{\lp}{\left(}
\newcommand{\rp}{\right)}
\newcommand{\lc}{\left[}
\newcommand{\rc}{\right]}

\newtheorem{theorem}{Theorem}[section]

\newtheorem{corollary}[theorem]{Corollary}

\newtheorem{definition}[theorem]{Definition}

\newtheorem{lemma}[theorem]{Lemma}

\newtheorem{proposition}[theorem]{Proposition}

\newtheorem{remark}[theorem]{Remark}

\begin{document}
\title{It\^o's formula for linear fractional PDEs}

\author{
  { Jorge A. Le\'{o}n}
   \thanks{Partially supported by the CONACyT grant 45684-F}    \\
  {\small\it Depto. de Control Autom\'{a}tico, CINVESTAV-IPN}   \\[-0.1cm]
  {\small\it Apartado Postal 14-740, 07000 M\'{e}xico, D.F., Mexico} \\[-0.1cm]
 {\small  {\tt jleon@ctrl.cinvestav.mx}}   
\and
  { Samy Tindel}              \\
  {\small\it Institut Elie Cartan, Universit\'e de Nancy 1} \\[-0.1cm]
  {\small\it BP 239 -- 54506 Vand\oe uvre-l\`es-Nancy, France}  \\[-0.1cm]
  {\small  {\tt tindel@iecn.u-nancy.fr }} \\[-0.1cm]  
  {\protect\makebox[5in]{\quad}}  
  \\
}

\maketitle

\begin{abstract}
In this paper we introduce a stochastic integral with respect to
the solution $X$ of the fractional heat equation on $\ou$, interpreted
as a divergence operator. This allows to use the
techniques of the Malliavin calculus in order to establish an 
It\^o-type formula for the process $X$. 
\end{abstract}

\vspace{1cm}


\vspace{0.3cm}


\section{Introduction}
In the last past years, a great amount of effort has been devoted
to a proper definition of  stochastic {\small PDE}s driven by a general noise.
For instance, the case of stochastic heat and wave equations in
$\R^n$ driven by a Brownian motion in time, with some mild conditions
on its spatial covariance, has been considered e.g. in 
\cite{Da,PZ,MS}, leading to some optimal results. More recently,
the case of {\small SPDE}s driven by a fractional Brownian motion has
been analyzed in \cite{Ca,TTV} in the linear case, or in \cite{GT,QT}
for the non-linear situation.

In this context, it seems natural to investigate the basic properties
(H\"olderianity, behavior of the density, invariant measures, 
numerical approximations, etc) of these objects. And indeed, in 
case of an equation driven by a Brownian motion, a lot of effort has 
been made in this direction (let us cite \cite{SMM,MS,GM} among others).
On the other hand, results concerning {\small SPDE}s driven by a 
fractional Brownian motion are rather scarce (see however \cite{NO}
for a result on {\small SPDE}s with irregular coefficients, and
\cite{SV} for a study of the H\"older regularity of solutions).

This article proposes then to go further into the study of processes
defined by fractional {\small PDE}s, and we will establish a 
It\^o-type formula for a random function $X$ on $\ott\times\ou$
defined as the solution to the heat equation with an additive
fractional noise. More specifically, we will consider $X$ as the
solution to the following equation:
\begin{equation}\label{intro:heat}
\partial_t X(t,x)=\Delta X(t,x) + B(dt,dx),
\quad (t,x)\in\ott\times\ou,\\
\end{equation}
with Dirichlet boundary conditions and null initial condition.
In equation (\ref{intro:heat}), the driving noise $B$ will be considered
as a fractional Brownian motion in time, with Hurst parameter
$H>1/2$, and as a white noise in space (notice that some more general
correlations in space could have been considered, as well as the
case $1/3<H<1/2$, but we have restrained ourselves to this simple situation
for sake of conciseness).

Then, for $X$ solution to (\ref{intro:heat}), $t\in\ott$, $x\in\ou$
and a $C_b^2$-function $f:\R\to\R$, we will prove that $f(X(t,x))$
can  be decomposed into:
\begin{equation}\label{intro:ito}
f(X(t,x))=f(0)
+\iot\iou \lp\ms_{t,x}f'(X)\rp(s,y) W(ds,dy)
+\frac12 \iot f''(X(s,x))K_{x}(ds),
\end{equation}
where in the last formula, $\ms_{t,x}$ is an operator based on the heat
kernel $G_t$ on $\ou$ and the covariance function of $B$, $W$ is 
a space-time white noise, and $K_x$ is the function defined on $\ott$ by:
$$
K_{x}(s)
=
H(2H-1) \int_0^s \int_0^s G_{2s-v_1-v_2}(x,x) |v_1-v_2|^{2H-2}
dv_1 dv_2.
$$
Notice also that, in (\ref{intro:ito}), the stochastic integral has to be
interpreted in the Skorohod sense (see Theorem \ref{teo:ito} for a 
precise statement).

As mentioned above, once the existence and uniqueness of the 
solution to (\ref{intro:heat}) is established, it certainly seems 
to be  a natural question to ask whether an It\^o-type
formula is available for the process we have produced. Furthermore,
this kind of result can also yield a better understanding of some properties
of the process itself, such as the distribution of hitting times, as
shown in \cite{DN}.
It is also worth mentioning at this point that
formula (\ref{intro:ito}) will be obtained thanks to some Gaussian
tools inspired by the case of the fractional Brownian motion itself. 
This is due to the fact that $X$ can be represented by the convolution
\begin{equation}\label{intro:convol}
X(t,x)=\iot\iou M_{t,s}(x,y) W(ds,dy)
\end{equation}
of a certain kernel $M$ on $\ott\times\ou$, defined at (\ref{def:mts}),
with respect to $W$. This kind of property
has already been exploited in \cite{GNT} for the case of the heat 
equation driven by a space-time white noise, but let us stress here
two differences with  respect to this latter reference: 
\begin{enumerate}
\item
On the one
hand, an important step of our computations will be to obtain the
representation (\ref{intro:convol}) itself (see Corollary \ref{cor:rep-x})
and to give some reasonable bounds on the kernel $M$ and its derivatives. 
\item
On the other
hand, the little gain in regularity we have in the current situation with 
respect to \cite{GNT} will allow us to obtain a formula for 
$t\mapsto f(X(t,x))$, while in the latter reference, we had to restrict 
ourselves to a change of variable formula for 
$$
t\mapsto \iou f(X(t,x)) \psi(x)\, dx,
$$
for a continuous function $\psi$.
\end{enumerate}

Let us say now a few words about the method we have used in order to get our 
result: as mentioned above, the first step in our approach consists in 
establishing the representation (\ref{intro:convol}). This representation,
together with the properties of the kernel $M$, suggest that the differential
of $X$ should be of the form
\begin{equation}\label{intro:diff-x}
X(dt,x)=
\lc \iot\iou \partial_t M_{t,s}(x,y) W(ds,dy) \rc dt.
\end{equation}
This formula is of course ill-defined, since 
$(s,y)\mapsto\partial_t M_{t,s}(x,y)$ is not a $L^2$-function on
$\ot\times\ou$, but it holds true for a regularization $M^\ep$
of $M$. We will then obtain easily an It\^o type formula for the process
$X^\ep$ corresponding to $M^\ep$, where the differential 
(\ref{intro:diff-x}) appears. Therefore, the  main step in our
 calculations will be to study the limit of the regularized It\^o formula
when $\ep\to 0$. Notice that this approach is quite different
(and from our point of view more intuitive) from the one adopted
in \cite{AMN1,GNT}, where the quantity $E[f(X(t,x)) I_n(\vp)]$
was evaluated for an arbitrary multiple integral $I_n(\vp)$
with respect to $W$.

Our paper is divided as follows: at Section 2, we will describe precisely
the noise and the equation under consideration, and we will give some
basic properties of the process $X$. Section 3 is devoted to the derivation
of our It\^o-type formula: at Section 3.1 we obtain the representation
(\ref{intro:convol}) for $X$, the regularized formula is given at 
Section 3.2, and eventually the limiting procedure is carried out  
at Sections 3.3
and 3.4. In the sequel of the paper, $c$ will designate a positive 
constant whose exact value can change from line to line.

\section{Preliminary definitions}
In this section we introduce the framework that will be used in this paper:
we will define precisely the noise which will be considered, then give
a brief review of some Malliavin calculus tools, and eventually
introduce the fractional heat equation.
\subsection{Noise under consideration}
Throughout the article, we will consider a complete probability space
$(\Omega,\cf,P)$ on which we
define a noise that will be a fractional Brownian motion with
Hurst parameter $H>1/2$ in time, and a Brownian motion in space.
More specifically, we define a zero mean Gaussian field 
$B=\{B(s,x):s\in[0,T],\ x\in[0,1]\}$ of the form
\begin{equation}\label{eq:n1.1}
B(t,x)=\int_0^t\int_0^x K_H(t,s)W(ds,dy).
\end{equation}  
Here $W$ is a two-parameter Wiener process and $K_H$ is the kernel
of the fractional Brownian motion (fBm) with Hurst parameter
$H\in (\frac12 ,1)$. Namely, for $0\le s\le t\le T$, we have
$$K_H(t,s)=C_Hs^{\frac12 -H}\int_s^t(u-s)^{H-\frac32}u^{H-\frac12}du,$$
where $C_H$ is a constant whose exact value is not important for our aim.
Observe that the standard theory of martingale measures introduced
in \cite{Wa} easily yields the existence of the integral (\ref{eq:n1.1}).

Note that it is natural to interpret the left-hand side of (\ref{eq:n1.1})
as the stochastic integral
\begin{equation}\label{eq:n2.0}
B(1_{[0,t]\times [0,x]}):=\int_0^t\int_0^x B(ds,dy).
\end{equation}
 The domain of this Wiener integral 
is then extended as follows: let ${\cal H}$ be the Hilbert
 space defined as the completion of the step functions with respect to
the inner product
\begin{eqnarray}\label{eq:n2.1}
\langle 1_{[0,s]},1_{[0,t]}\rangle_{\cal H}&=&
\langle K_H(t,\cdot),K_H(s,\cdot)\rangle_{L^2([0,T])}\nonumber\\
&=&H(2H-1)\int_0^t\int_0^s|u-r|^{2H-2}dudr.
\end{eqnarray}
Thus, by Al\`os and Nualart \cite{Alda}, the kernel $K_H$ allows to 
construct an isometry $K_{H,T}^*$ from ${\cal H}\times L^2([0,1])$
(denoted by  ${\cal H}_T$ for short) into $L^2([0,T]\times [0,1])$
such that, for $0\le s <t\le T$,
\begin{eqnarray*}
\lp K^*_{H,T}1_{[0,t]\times[0,x]}\rp(s,y)
&&=K_H(t,s)1_{[0,x]}(y)\\
&&=1_{[0,x]}(y)\int_s^T1_{[0,t]}(r)\partial_r K_H
(r,s)dr.
\end{eqnarray*}
Therefore the Wiener integral (\ref{eq:n2.0}) can be extended into
an isometry $\varphi \mapsto B(\varphi)$ from ${\cal H}_T$ into a
subspace of $L^2(\Omega)$ so that, for any $\vp\in\ch_T$,
\begin{equation}\label{eq:n3.1}
B(\varphi)=\int_0^T\iou (K^*_{H,T}\varphi)(s,y)W(ds,dy).
\end{equation}
Then, for two elements $\vp$ and $\psi$ of $\ch_T$, the covariance
between $B(\vp)$ and $B(\psi)$ is given by
\begin{equation}\label{cov:str}
E\lc B(\vp) B(\psi)  \rc
=H(2H-1)
\iott\iott\iou \vp(s,y) |s-r|^{2H-2} \psi(r,y)\, ds dr dy.
\end{equation}

Notice that
an element of ${\cal H}_T$ could possibly not be a function. 
Hence, as the in fBm case, we will deal with 
the Banach space $|{\cal H}_T|$ of all the measurable functions
$\varphi:[0,T]\times[0,1]\rightarrow \R$ such that
\begin{eqnarray}\label{eq:n3.2}
||\varphi||_{|{\cal H}_T|}&=&H(2H-1)\int_0^T\int_0^T\iou |\varphi
(r,y)||u-r|^{2H-2}|\varphi(u,y)|dydudr\nonumber \\
&=&\iou\int_0^T\lp\int_s^T|\varphi(r,y)|\partial_r K_H
(r,s)dr\rp^2dsdy<\infty\nonumber . 
\end{eqnarray}
It is then easy to see that $L^2([0,T]\times[0,1])\subset
|{\cal H}_T|\subset {\cal H}_T$.

\subsection{Malliavin calculus tools}
The goal of this section is to recall the basic definitions of 
 the Malliavin calculus  which will allow us
to define the divergence operator with respect to $W$. For
a more detailed presentation, we recommend  Nualart \cite{Nul}.

Let ${\cal S}$ be the family of all smooth functionals of the form
$$
F=f(W(s_1,y_1),\ldots , W(s_n,y_n)),
\quad\mbox{ with }\quad
(s_i,y_i)\in\ott\times\ou,
$$
where $f\in C_b^{\infty}(\R^n)$ (i.e., $f$ and all its partial derivatives
are bounded). The derivative of this kind of smooth functional is the
$L^2([0,T]\times[0,1])$-valued random variable
$$DF=\sum_{i=1}^n \frac{\partial f}{\partial x_i}
(W(s_1,y_1),\ldots , W(s_n,y_n))1_{[0,s_i]\times[0,y_i]}.
$$
It is then well-known that $D$ is a closeable operator from $L^2(\Omega)$ into
$L^2(\Omega\times[0,T]\times[0,1])$. Henceforth, to simplify the
notation, we also denote its closed extension by $D$. Consequently
$D$ has an adjoint $\delta$, which is also a closed operator, 
characterized via the duality relation
$$E\lp F\delta(u)\rp=E\lp \langle DF,u\rangle_{L^2([0,T]\times[0,1])}
\rp,$$
with $F\in{\cal S}$ and $u\in\dom(\delta)
\subset L^2(\Omega\times[0,T]\times[0,1])$.
The operator $\delta$ has been considered as a stochastic integral 
because it is an extension of the It\^o integral with respect to
$W$ that allows us to integrate anticipating processes
(see, for instance, \cite{Nul}). According to this fact, we will
sometimes use the notational convention
$$\delta(u)=\int_0^T\iou u_{s,y}W(ds,dy).$$

Notice that the operator $\delta$ (or Skorohod integral) 
has the following property:
Suppose that $F$ is a random variable in $\dom(D)$ and that $u$ is Skorohod
integrable (i.e., $u\in\dom(\delta)$), such that 
$E(F^2\int_0^T\iou (u(s,y))^2dyds)<\infty$. Then 
\begin{equation}\label{eq:n5.0}
\int_0^T\iou Fu(s,y)W(ds,dy)=F\int_0^T\iou
u(s,y)W(ds,dy)-\int_0^T\iou(D_{s,y}F)u(s,y)dyds,
\end{equation}
in the sense that $(Fu)\in\dom(\delta)$ if and only if
the right-hand side is in $L^2(\Omega)$.

\subsection{Heat equation}
This paper is concerned with the solution $X$ to the following stochastic
heat equation on $\ou$, 
with Dirichlet boundary conditions and null initial condition:
\begin{equation}\label{eq:heat}
\begin{cases}
\partial_t X(t,x)=\Delta X(t,x) + B(dt,dx),& (t,x)\in\ott\times\ou\\
X(0,x)=0,\quad X(t,0)=X(t,1)=0.&
\end{cases}
\end{equation}
It is well-known (see \cite{TTV}) that equation (\ref{eq:heat}) has a unique
solution, which is given explicitly by
\begin{equation}\label{sol:heateq}
X(t,x)= \iot\iou G_{t-s}(x,y) B(ds,dy),
\end{equation}
where
\begin{equation}\label{eq:n5.1}
G_t(x,y)=\frac{1}{\sqrt{4\pi t}}\sum_{n=-\infty}^{\infty}
\lc \exp\lp -\frac{(y-x-2n)^2}{4t}\rp-\exp\lp -\frac{(y+x-2n)^2}{4t}\rp\rc
\end{equation}
stands for the Dirichlet heat kernel on $\ou$
with Dirichlet boundary conditions.
Let us recall here some elementary but useful identities
for the heat kernel $G$:
\begin{lemma}\label{heat-ker}
The following relations hold true for the heat kernel $G$ given by
(\ref{eq:n5.1}):
$$
\iou G_t(x,y) dy =1,\quad
G_t(x,y)\le
\frac{c_1}{t^{1/2}}\exp\lp -\frac{c_2(x-y)^2}{t} \rp,
$$
and
$$
|\partial_t G_t(x,y)|\le
\frac{c_3}{t^{3/2}}\exp\lp -\frac{c_4(x-y)^2}{t} \rp,
$$ for some positive constants $c_1,\ c_2,\ c_3$ and $c_4$.
Furthermore, $G$ can be decomposed into
\begin{equation}\label{dcp:ker}
G_t(x,y)=G_{1,t}(x,y)+ R_t(x,y),
\end{equation}
where
$$
G_{1,t}(x,y)
=\frac{1}{\sqrt{4\pi t}}
\lc \exp\lp -\frac{(y-x)^2}{4t}\rp-\exp\lp -\frac{(y+x)^2}{4t}\rp-
\exp\lp -\frac{(y+x-2)^2}{4t}\rp\rc,
$$
and
$R_t(x,y)$ is a smooth bounded function on  $\ott\times\ou^2$.
\end{lemma}

Let us recall now some basic properties of the process $X$ defined 
by (\ref{eq:heat}) and (\ref{sol:heateq}), starting with its integrability.
\begin{lemma}
The process defined on $\ott\times\ou$ by (\ref{sol:heateq})
satisfies
$$
\sup_{t\in\ott,x\in\ou} E\lc |X(t,x)|^2  \rc <\infty.
$$
\end{lemma}

\begin{proof}
We have, according to (\ref{cov:str})
 and Lemma \ref{heat-ker}, that
\begin{eqnarray*}
E\lc |X(t,x)|^2  \rc&=&
c_H \iotd \frac{ds du}{|s-u|^{2-2H}}
\iou G_{t-s}(x,y) G_{t-u}(x,y) dy\\
&\le& c
\iotd \frac{ds du}{(t-s)^{1/2}|s-u|^{2-2H}}
\iou  G_{t-u}(x,y) dy\\
&=&
c \iotd \frac{ds du}{(t-s)^{1/2}|s-u|^{2-2H}},
\end{eqnarray*}
and the last integral is finite by elementary arguments.
\end{proof}
\vspace{0.3cm}

One can go further in the study of $X$, and show the following
regularity result (see also \cite{SV}):
\begin{proposition}\label{lem:n1.3}
Let $X$ be the solution to (\ref{eq:heat}). Then, for $t_1,t_2\in\ott$
and $x\in\ou$, we have
$$
E\lc |X(t_2,x)-X(t_1,x)|^2  \rc
\le
c |t_2-t_1|^{2\ga},
$$
for any $\ga<H-1/4$. In particular, for any $T>0$ and $x\in\ou$, the
function $t\in\ott\mapsto X(t,x)$ is $\ga$-H\"older continuous for
any $\ga<H-1/4$.
\end{proposition}

\begin{proof}
Assume $t_1<t_2$. We then have
$$
X(t_2,x)-X(t_1,x)
=
A(t_1,t_2,x) + B(t_1,t_2,x),
$$
with
$$A(t_1,t_2,x)=
\int_0^{t_1}\iou \lc G_{t_2-s}(x,y) - G_{t_1-s}(x,y)\rc B(ds,dy)$$
and
$$B(t_1,t_2,x)=
\int_{t_1}^{t_2}\iou G_{t_2-s}(x,y)  B(ds,dy).
$$
Hence
\begin{equation}\label{eq:n8.1}
E\lc |X(t_2,x)-X(t_1,x)|^2  \rc
\le
2 \lp E\lc A^2(t_1,t_2,x)  \rc + E\lc B^2(t_1,t_2,x)  \rc  \rp.
\end{equation}

We first note that (\ref{cov:str}) and Lemma \ref{heat-ker}
imply
\begin{eqnarray}\label{eq:n8.2}
\lefteqn{E\lc B^2(t_1,t_2,x)\rc}\nonumber\\
&=&c_H\int_{t_1}^{t_2}\int_{t_1}^{t_2}duds|s-u|^{2H-2}\iou
G_{t_2-s}(x,y)G_{t_2-u}(x,y)dy\nonumber\\
&&\le c\int_{t_1}^{t_2}ds(t_2-s)^{-1/2}\int_{t_1}^{t_2}|s-u|^{2H-2}du
\nonumber\\
&&\le c(t_2-t_1)^{2H-\frac12}.
\end{eqnarray}

Now we will concentrate  on the estimate on
$E[ A^2(t_1,t_2,x) ]$. By (\ref{cov:str}), we have
\begin{equation}\label{eq:n8.3}
E\lc A^2(t_1,t_2,x)  \rc
=c_H
\int_0^{t_1}\int_0^{t_1}
\frac{du ds}{|s-u|^{2-2H}} C_x(s,u),
\end{equation}
with $C_x(s,u)$ defined by
$$
C_x(s,u)=
\iou \lc G_{t_2-s}(x,y) - G_{t_1-s}(x,y)\rc
\lc G_{t_2-u}(x,y) - G_{t_1-u}(x,y)\rc dy.
$$
Thus, invoking Lemma \ref{heat-ker}, we obtain that, for a given
$\al<1/2$,
$$
C_x(s,u)
\le c
\frac{(t_2-t_1)^{2\al}}{(t_1-u)^{3\al/2}(t_1-s)^{3\al/2}}
D_x(s,u),
$$
where
$$
D_x(s,u)
=\iou
|G_{t_2-s}(x,y) -G_{t_1-s}(x,y)|^{1-\al}
| G_{t_2-u}(x,y) - G_{t_1-u}(x,y)|^{1-\al}
 dy.
$$
It is then easily seen that $D_x(s,u)$ can be bounded by a sum of
terms of the form
$$
F_x(s,u)=
\iou G_{\si-s}^{1-\al}(x,y) G_{\tau-u}^{1-\al}(x,y) dy,
$$
with $\si,\tau\in\{ t_1,t_2 \}$. This latter expression can be bounded
in the following way:
\begin{eqnarray*}
F_x(s,u)
&\le&
\lp \iou G_{\si-s}^{2(1-\al)}(x,y) dy \rp^{1/2}
\lp \iou G_{\tau-u}^{2(1-\al)}(x,y) dy \rp^{1/2}\\
&\le&
\frac{c}{(t_1-s)^{1/4-\al/2}(t_1-u)^{1/4-\al/2}}.
\end{eqnarray*}
We have thus obtained that
$$
E\lc A^2(t_1,t_2,x)  \rc
\le
c(t_2-t_1)^{2\al}
\int_0^{t_1}\int_0^{t_1}
\frac{du ds}{|s-u|^{2-2H}(t_1-s)^{1/4+\al}(t_1-u)^{1/4+\al}}.
$$
Now thanks the change of variable $v=\frac{u-s}{t_1-s}$, 
 the latter integral is finite
whenever $\al<H-1/4$, which, together with (\ref{eq:n8.1}) and
(\ref{eq:n8.2}), ends the proof.
\end{proof}

\section{It\^o's formula for the heat equation}
Let us turn to the main aim of this paper, namely the It\^o-type
formula for the process $X$ introduced
in (\ref{sol:heateq}). The strategy of our computations
can be briefly outlined as follows: first we will try to represent
$X$ as a convolution of a certain kernel $M$ with respect to $W$,
with reasonable bounds on $M$. Then we will be able to establish our
It\^o's formula for a smoothed version of $X$, involving a regularized
kernel $\me$ for $\ep>0$, by applying the usual It\^o formula.
Our main task will then be to study the limit of the quantities
we will obtain as $\ep\to 0$.

\subsection{Differential of $X$}

Before getting a suitable expression for the differential of
$X$, let us see how to represent this process as a convolution
with respect to $W$.

\subsubsection{Representation of $X$}\label{sec:rep-x}

The expressions (\ref{eq:n2.1}) and (\ref{eq:n3.1}) lead to the following
result (see \cite{Alda}).
\begin{lemma}\label{lem:n2.1}
Let  $\vp$ be a function in $|{\cal H}_T|$. Then 
$$
\iot\iou \vp(s,y) B(ds,dy)
=
\iot\iou [K^*_{H,T}1_{[0,t]}\vp](u,y) W(du,dy),
$$
with
$$
[K^*_{H,T}1_{[0,t]}\vp](u,y)
=1_{[0,t]}(u)
\int_u^t \vp(r,y) \partial_r \kh(r,u) dr.
$$
\end{lemma}
\begin{remark} \label{re:heur}
This result could also have been obtained by
some heuristic arguments. Indeed, a formal
 way to write (\ref{eq:n1.1})
is to say that, for $t>0$ and $y\in\ou$, the differential $B(t,dy)$ is
defined as
$$
B(t,dy)=\iot \kh(t,s) W(ds,dy).
$$
Thus, if we differentiate formally this expression in time,
since $\kh(t,t)=0$, we obtain
$$
\partial_t B(t,dy)
=
\lc \iot \partial_t \kh(t,s) W(ds,dy) \rc dt.
$$
Since
$\partial_t \kh(t,s)$ is not a $L^2$-function, the
last equality has to be interpreted
in the following way: if $\vp$ is a deterministic function,
then
\begin{eqnarray*}
\iot\iou \vp(s,y) B(ds,dy)
&=&
\iot\iou \vp(s,y)
\lc \int_0^s \partial_s \kh(s,u) W(du,dy) \rc ds\\
&=&\iot\iou W(du,dy)
\lc \int_u^t \vp(s,y) \partial_s K(s,u) ds \rc,
\end{eqnarray*}
which recovers the result of Lemma \ref{lem:n2.1}.
\end{remark}
We can now easily get the announced representation for $X$:
\begin{corollary}\label{cor:rep-x}
The solution $X$ to (\ref{eq:heat}) can be written as
\begin{equation}\label{rep:x}
X(t,x)=\iot\iou M_{t,s}(x,y) W(ds,dy),
\end{equation}
with
\begin{equation}\label{def:mts}
M_{t,s}(x,y)= \ist G_{t-u}(x,y) \partial_u \kh(u,s) du.
\end{equation}
\end{corollary}
\begin{proof} The result is an immediate consequence of the proof of
Proposition \ref{lem:n1.3} and Lemma \ref{lem:n2.1}.\end{proof}

\subsubsection{Some bounds on $M$}

The kernel $M$ will be algebraically useful in order to obtain
our It\^o's formula, and we will proceed to show now that it behaves
similarly to the heat kernel $G$. To do so, let us first state the following
technical lemma:
\begin{lemma}\label{est:der-f}
Let $f$ be defined on $0<r<t\le T$ by
$$
f(r,t)=
\int_r^t (t-u)^{-1/2}(u-r)^{-\al}\exp\lp -\frac{\kappa x^2}{t-u} \rp du,
$$
for a constant $\ka>0$, $x\in[0,2]$ and $\alpha\in(0,1)$.
Then, there exist some constants $c_1,c_2,c_3,c_4>0$ such that
\begin{equation}\label{majdtf}
f(r,t)\le c_1  (t-r)^{-(\al-1/2)}
\exp\lp -\frac{c_2 x^2}{t-r} \rp\end{equation}
and
\begin{equation}\label{eq:n12.1}
\partial_tf(r,t)\le c_3  (t-r)^{-(\al+1/2)}
\exp\lp -\frac{c_4 x^2}{t-r} \rp .
\end{equation}
\end{lemma}

\begin{proof}
Recall that, in the remainder of the paper, $\ka$ stands for a positive 
constant which can change from line to line.
Notice also that (\ref{majdtf}) is easy to see due to
$$f(r,t)\le \exp\lp -\frac{\ka x^2}{t-r}\rp\int_r^t 
(t-u)^{-1/2}(u-r)^{-\alpha}dr.$$
Now we will concentrate on (\ref{eq:n12.1}):
let us perform the change of variable $v=\frac{u-r}{t-r}$. This yields
$$
f(r,t)=
(t-r)^{-(\al-1/2)}\int_0^{1}
(1-v)^{-1/2}v^{-\al}
\exp\lp -\frac{\kappa x^2}{(1-v)(t-r)} \rp dv,
$$
and thus
$$
\partial_t f(r,t)=g_1(r,t)+g_2(r,t),
$$
with
$$
g_1(r,t)=\kappa x^2
(t-r)^{-(\al+3/2)}\int_0^{1}
(1-v)^{-3/2}v^{-\al} 
\exp\lp -\frac{\ka x^2}{(1-v)(t-r)} \rp dv $$ and
$$g_2(r,t)=
\lp \frac12 -\al\rp
(t-r)^{-(\al+1/2)}\int_0^{1}
(1-v)^{-1/2}v^{-\al}
\exp\lp -\frac{\ka x^2}{(1-v)(t-r)} \rp dv.
$$
Therefore, thanks to the fact that $u\mapsto ue^{-u}$ is a bounded
function on $\R_+$, we have
\begin{eqnarray*}
g_1(r,t)
&\le&
c (t-r)^{-(\al+1/2)}\int_0^{1}
(1-v)^{-1/2}v^{-\al} 
\exp\lp -\frac{\ka x^2}{2(1-v)(t-r)} \rp dv
\\
&\le&c(t-r)^{-(\al+1/2)}\exp\lp -\frac{\ka x^2}{2(t-r)}\rp
\iou (1-v)^{-1/2}v^{-\al}dv,
\end{eqnarray*}
which is an estimate of the form (\ref{majdtf}).
Finally, it is easy to see that
$$g_2(r,t)
\le c(t-r)^{-(\al+1/2)}\exp\lp -\frac{\ka x^2}{2(t-r)}\rp
\iou (1-v)^{-1/2}v^{-\al}dv,
$$
which completes the proof.

\end{proof}

\vspace{0.3cm}

We are now ready to prove our bounds on $M$:
\begin{proposition}\label{prop:n2.4}
Let $M$ be the kernel defined at (\ref{def:mts}). Then, for some strictly
positive constants $c_5,c_6,c_7,c_8>0$,
we have
\begin{eqnarray*}
\lefteqn{M_{t,s}(x,y)}\\
 &\le&
c_5 (t-s)^{-(1-H)} \lp\frac{t}{s}\rp^{H-1/2}
\lc \exp\lp -\frac{c_6 (x-y)^2}{t-s} \rp+
\exp\lp -\frac{c_6 (x+y-2)^2}{t-s} \rp\rc
\end{eqnarray*}
and
\begin{eqnarray*}
\lefteqn{|\partial_t M_{t,s}(x,y)|}\\
 &\le&
c_7 (t-s)^{-(2-H)} \lp\frac{t}{s}\rp^{H-1/2}
\lc\exp\lp -\frac{c_8 (x-y)^2}{t-s} \rp
+\exp\lp -\frac{c_8 (x+y-2)^2}{t-s} \rp\rc.
\end{eqnarray*}
\end{proposition}

\begin{proof}
 First of all, we will use the decomposition
(\ref{dcp:ker}), which allows to write
$$
M_{t,s}(x,y)
=\ist G_{1,t-u}(x,y) \partial_u \kh(u,s) du
+
\ist R_{t-u}(x,y) \partial_u \kh(u,s) du.
$$
Now the result is an immediate consequence of Lemma \ref{est:der-f}
applied to $\alpha<\frac32-H$, the only difference
being the presence of the term $(u/s)^{H-1/2}$, which
can be bounded by $(t/s)^{H-1/2}$ each time it appears.
This yields the desired result.

\end{proof}

\subsubsection{Differential of $X$}\label{sec:dx}

With the representation (\ref{rep:x}) in hand, we can now follow
the heuristic steps in Remark \ref{re:heur} in order
to get a reasonable definition of the differential of
$X$ in time. That is, we can write formally that
$$
X(dt,x)=\lc \iot\iou \partial_t M_{t,s}(x,y) W(ds,dy) \rc dt ,
$$
which means that if $\vp:[0,T]\times [0,1]\rightarrow \R$ is a smooth enough
function, we have
\begin{eqnarray*}
\iott\vp(t,x) X(dt,x)
&=&
\iott \vp(t,x)
\lc \iot\iou \partial_t M_{t,s}(x,y) W(ds,dy) \rc dt \\
&=&
\iott\iou W(ds,dy)
\lc \int_s^T  \vp(t,x)\partial_t M_{t,s}(x,y) dt \rc.
\end{eqnarray*}
Note that this expression may not be convenient because it does not
take advantage of the continuity of $\vp$. But, by 
Proposition \ref{prop:n2.4}, we can write
$$\int_s^T  \vp(t,x)\partial_t M_{t,s}(x,y) dt=
\int_s^T (\vp(t,x)-\vp(s,x))\partial_t M_{t,s}(x,y) dt+\vp(s,x)
M_{T,s}(x,y).$$
Here again, we can formalize these heuristic considerations
into the following:
\begin{definition}
Let $\vp:\Omega\times [0,T]\times\ou\rightarrow\R$ be a measurable process. 
We say
that $\vp$ is integrable with respect to $X$ if the mapping 
\begin{equation}\label{def:mtx-phi}
(s,y)\mapsto [\ms_{T,x}\vp](s,y)
:=
\int_s^T \left(\vp(t,x)-\vp(s,x)\right) \partial_t M_{t,s}(x,y) dt
+\vp(s,x)M_{T,s}(x,y)
\end{equation}
belongs
 to $\dom(\delta)$, 
for almost all $x\in\ou$. In this case we set
$$
\iott \vp(t,x) X(dt,x)
=
\iott\iou [\ms_{T,x}\vp](s,y) W(ds,dy).
$$
\end{definition}
\begin{remark}\label{rmk:interp-div}
Just like in the case of the fractional  Brownian motion \cite{AMN1}
or of the heat equation driven by the space-time white noise \cite{GNT},
one can show that $\iott \vp(t,x) X(dt,x)$ can be interpreted as a 
divergence operator for the Wiener space defined by $X$.
\end{remark}
\begin{remark}\label{rmk:holder}
It is easy to see that Proposition \ref{prop:n2.4} implies that
$\vp: [0,T]\rightarrow\R$ is integrable with respect to $X$ if it is 
$\beta$-H\"older continuous in time with $\beta>1-H$.
\end{remark} 

\subsection{Regularized version of It\^o's formula}\label{sec:regv}

The representation (\ref{rep:x}) of $X$ also allows us to define
a natural regularized version $\xe$ of $X$, depending
on a parameter $\ep>0$, such that $t\mapsto \xe(t,x)$ will be
a semi-martingale. Indeed, set, for $\ep>0$,
$$
\me_{t,s}(x,y)
=
\ist G_{t-u+\ep}(x,y) \partial_u\kh(u+\ep,s) du,
$$
and
\begin{equation}\label{def:xep}
\xe(t,x)=\iot\iou \me_{t,s}(x,y) W(ds,dy).
\end{equation}
We will also need a regularized operator $\mes_{t,x}$
(see (\ref{def:mtx-phi})), defined
naturally by
$$
\lc \mes_{t,x}\vp\rc(s,y)
=
\int_s^t (\vp(r,x)-\vp(s,x)) \partial_r \me_{r,s}(x,y) dr 
+\vp(s,x)\me_{t,s}(x,y).
$$
Our strategy in order to get an It\^o type formula for $X$
will then be the following:
\begin{enumerate}
\item
Apply the usual It\^o formula to the semi-martingale $t\mapsto \xe(t,x)$.
\item
Rearrange terms in order to get an expression in terms
of the operator $\mes_{t,x}$.
\item
Study the limit of the different terms obtained through Steps
1 and 2, as $\ep\to 0$.
\end{enumerate}
The current section will be devoted to the elaboration of Steps
1 and 2. 
\begin{lemma}\label{lem:n2.6'} 
Let $\ep>0$. Then, the process $t\mapsto\xe(t,x)$ has bounded variations on 
$[0,T]$, for all $x\in\ou$.\end{lemma}
\begin{proof}
The Fubini theorem for $W$ and the semigroup property of $G$ imply
$$
\xe(t,x)=\int_0^t\iou G_{t-u+\frac{\ep}{2}}(x,z)
\lp\int_0^u\iou G_{\ep/2}(z,y)\partial_u\kh(u+\ep,s)W(ds,dy)\rp dzdu,
$$
and notice that this integral is well-defined due to Kolmogorov's
continuity theorem.
Therefore, since $t\mapsto G_{t-u+\ep/2}(x,z)$ is also a $C^1$-function
on $[u,T]$, we obtain that 
$\xe$ is differentiable  with respect to $t\in[0,T]$, and
\begin{eqnarray*}
\lefteqn{\partial_t\xe(t,x)} \\
&=&\int_0^t\iou \partial_t G_{t-u+\frac{\ep}{2}}(x,z)
\lp\int_0^u\iou G_{\ep/2}(z,y)\partial_u\kh(u+\ep,s)W(ds,dy)\rp dzdu\\
&&+\iou G_{\ep/2}(x,z)
\lp\int_0^t\iou G_{\ep/2}(z,y)\partial_t\kh(t+\ep,s)W(ds,dy)\rp dz,
\end{eqnarray*}
which is a continuous process on $\ott\times\ou$,
invoking Kolmogorov's continuity theorem again in a standard manner.

\end{proof}

\noindent
An immediate consequence of the previous lemma is the following:
\begin{corollary}\label{cor:n2.6''}
Let $t\in[0,T]$, $x\in [0,1]$ and $\ep>0$. Then,
\begin{eqnarray*}
\partial_t\xe(t,x)
&=&\int_0^t\iou \lp\int_s^t\partial_t G_{t-u+\ep}(x,y)
\partial_u\kh(u+\ep,s)du\rp W(ds,dy)\\
&&+\int_0^t\iou G_{\ep}(x,y)
\partial_t\kh(t+\ep,s)W(ds,dy)\\
&=&\int_0^t\iou\partial_t M_{t,s}^{\ep}(x,y)W(ds,dy).
\end{eqnarray*}
\end{corollary}
\begin{proof}
The result follows from Fubini's theorem  for $W$ and from the
semigroup property of $G$.

\end{proof} 

Now we are ready to establish our regularized It\^o's formula 
in order to carry out Steps 1 and 2 of this section. 
\begin{proposition}\label{pro:n2.6}
Let $f$ be a regular function in $C_b^2(\R)$, $\ep>0$,
and $\xe$ the process defined by (\ref{def:xep}). Then, for
$t\in[0,T]$ and $ x\in\ou$, $\mes_{t,x}f'(\xe)$ belongs to 
$\dom(\delta)$ and
$$
f(\xe(t,x))=f(0)+\mathbf{A_{1,\ep}}(t,x)+\mathbf{A_{2,\ep}}(t,x),
$$
where
$$
\mathbf{A_{1,\ep}}(t,x)=
\iot\iou \lp\mes_{t,x}f'(\xe)\rp(s,y) W(ds,dy)
$$
is defined as a Skorohod integral, and
$$
\mathbf{A_{2,\ep}}(t,x)=
 \iot f^{''}(\xe(s,x))K_{\ep,x}(ds),
$$
with
\begin{eqnarray}\label{def:kex}
K_{\ep,x}(s)
&=&
\int_0^s dv_2\int_0^{v_2}dv_1 G_{2(s+\ep)-v_1-v_2}(x,x)
\Big\{ H(2H-1)|v_1-v_2|^{2H-2}  \nonumber\\
&&\qquad
-\partial^2_{v_1,v_2}\lp \int_{v_1}^{v_1+\ep}\kh(v_1+\ep,u)\kh(v_2+\ep,u)
du\rp\nonumber\\
&&\qquad - \partial_{v_2}\lp \kh(v_1+\ep,v_1)\kh(v_2+\ep,v_1)\rp
\Big\}.
\end{eqnarray}
\end{proposition}

\begin{proof}
By Corollary \ref{cor:n2.6''}, we  are able to apply 
the classical change of variable
formula to obtain
\begin{equation}\label{cv:class}
f(\xe(t,x))=f(0)+
\int_0^t f'(\xe(s,x))\lc \int_0^s\iou \partial_s \me_{s,u}(x,y) W(du,dy)
 \rc ds.
\end{equation}
Moreover, the   derivative of $f'(\xe(s,x))$ in the 
Malliavin calculus sense  is given by
$$
D_{v,z}[f'(\xe(s,x))]
=
\me_{s,v}(x,z)f''(\xe(s,x)) \1_{\{ v\le s\}}.
$$
Since the last quantity is bounded by $c_\ep v^{\frac12 -H}$ for
$\ep>0$, then invoking formula (\ref{eq:n5.0}) for the Skorohod integral, 
we get
\begin{align}\label{div:appli}
&f'(\xe(s,x)) \int_0^s\iou \partial_s \me_{s,u}(x,y) W(du,dy)\nonumber\\
=&
\int_0^s\iou f'(\xe(s,x)) \partial_s \me_{s,u}(x,y) W(du,dy)\nonumber\\
&\qquad\qquad\qquad +f''(\xe(s,x))
\int_0^s\iou \lp\partial_s \me_{s,u}(x,y)\rp \me_{s,u}(x,y) dudy.
\end{align}
Denote for the moment the quantity
$\int_0^s\iou \lp\partial_s \me_{s,u}(x,y)\rp \me_{s,u}(x,y) du dy$
by $h_x(s)$. Then, combining (\ref{cv:class}) and (\ref{div:appli}),
proceeding as the beginning of Section \ref{sec:dx},
and applying 
Fubini's theorem for the Skorohod integral,
we have 
\begin{equation}\label{eq:first-ito}
f(\xe(t,x))=f(0)+\mathbf{A_{1,\ep}}(t,x)+
\iot f''(\xe(s,x)) h_x(s) ds.
\end{equation}
We can find now a simpler expression for $h_x(s)$. Indeed, since
$\me_{s,s}(x,y)=0$, it is easily checked that
\begin{equation}\label{eq:hxs}
h_x(s)=
\frac12 \partial_s \lc \int_0^s\iou  \lp \me_{s,u}(x,y) \rp^2 du dy\rc.
\end{equation}
Furthermore, the semigroup property for $G$ yields
\begin{align*}
& \int_0^s\iou  \lp \me_{s,u}(x,y) \rp^2 du dy\nonumber\\
&= \int_0^s du \int_u^s dv_1 \int_u^s dv_2 \iou dy\,
G_{s+\ep-v_1}(x,y) G_{s+\ep-v_2}(x,y)
\partial_{v_1}\kh(v_1+\ep,u)\nonumber\\
&\hspace{11cm}\cdot \partial_{v_2}\kh(v_2+\ep,u),\nonumber\\
\end{align*}
and this last expression is equal to
\begin{align}\label{inter1:hxs}
&2\int_0^s du \int_u^s dv_1 \int_{v_1}^s dv_2 
G_{2(s+\ep)-v_1-v_2}(x,x)
\lp\partial_{v_1}\kh(v_1+\ep,u)\rp \partial_{v_2}\kh(v_2+\ep,u)\nonumber\\
&=2\int_0^s dv_2 \int_0^{v_2} dv_1
G_{2(s+\ep)-v_1-v_2}(x,x)
\lp \int_0^{v_1} \lp\partial_{v_1}\kh(v_1+\ep,u)\rp 
\partial_{v_2}\kh(v_2+\ep,u) du \rp.
\end{align}
But
\begin{eqnarray}\label{inter2:hxs}
\lefteqn{\int_0^{v_1} \lp\partial_{v_1}\kh(v_1+\ep,u)\rp
 \partial_{v_2}\kh(v_2+\ep,u) du}\nonumber\\
&=&
\partial_{v_2}\partial_{v_1} \lc 
\int_0^{v_1} \kh(v_1+\ep,u) \kh(v_2+\ep,u) du\rc\nonumber\\
&&-\partial_{v_2}\lc\kh(v_1+\ep,v_1) \kh(v_2+\ep,v_1)\rc\nonumber\\
&=& 
H(2H-1) |v_1-v_2|^{2H-2}
-\partial_{v_2}\partial_{v_1}\lc\int_{v_1}^{v_1+\ep}
\kh(v_1+\ep,u) \kh(v_2+\ep,u) du\rc\nonumber\\
&&-\partial_{v_2}\lc\kh(v_1+\ep,v_1) \kh(v_2+\ep,v_1)\rc .
\end{eqnarray}
By putting together (\ref{inter1:hxs}) and (\ref{inter2:hxs}), we have
thus obtained that
$$
\frac12 \int_0^s\iou  \lp \me_{s,u}(x,y) \rp^2 du dy
=
K_{\ep,x}(s),
$$
where $K_{\ep,x}(s)$ is defined at (\ref{def:kex}). By plugging this
equality into (\ref{eq:first-ito}) and (\ref{eq:hxs}), the proof is now
complete.

\end{proof}

\subsection{It\^o's formula}

We are now ready to perform the limiting procedure which will
allow to go from Proposition \ref{pro:n2.6} to the announced It\^o
formula. To this end we will need the following technical result,
which states that the modulus of continuity of $t\mapsto X^\ep(t,x)$
can be bounded from below by any $\nu<H-1/4$, independently of $\ep$.
\begin{proposition}\label{est:xep}
Let $\xe$ be given by (\ref{def:xep}). Then for $t_1,t_2\in[0,T]$ and
$x\in\ou$, there is a positive constant $c$ (independent of $\ep$) such that
$$E\lp|\xe(t_2,x)-\xe(t_1,x)|^2\rp\le c|t_2-t_1|^{2\nu},$$
for any $\nu<H-\frac14$.
\end{proposition}

\begin{proof} Suppose that $t_1<t_2$. Then
\begin{eqnarray}\label{eq:n21.1}
E\lp|\xe(t_2,x)-\xe(t_1,x)|^2\rp&\le&2
\int_0^{t_1}\iou\lp\me_{t_2,s}(x,y)- \me_{t_1,s}(x,y)\rp^2dyds\nonumber\\
&&+2\int_{t_1}^{t_2}\iou\lp\me_{t_2,s}(x,y)\rp^2dyds.
\end{eqnarray}
Now using the fact that $\partial_u\kh(u,s)>0$, 
we have
\begin{eqnarray}\label{eq:n21.2}
\lefteqn{\int_{t_1}^{t_2}\iou\lp\me_{t_2,s}(x,y)\rp^2dyds}\nonumber\\
&=&\int_{t_1}^{t_2}\iou\lp\int_{s+\ep}^{t_2+\ep}
G_{t_2+2\ep-u}(x,y)\partial_u\kh(u,s)du\rp^2dyds\nonumber\\
&&\le \int_{0}^{t_2+\ep}\iou\lp\int_{s}^{t_2+\ep}1_{[t_1+\ep,
t_2+\ep]}(u)
G_{t_2+2\ep-u}(x,y)\partial_u\kh(u,s)du\rp^2dyds\nonumber\\
&=& H(2H-1)
\int_{t_1+\ep}^{t_2+\ep}\int_{t_1+\ep}^{t_2+\ep}\iou |u-v|^{2H-2}
G_{t_2+2\ep-u}(x,y)G_{t_2+2\ep-v}(x,y)dydudv\nonumber\\
&&\le c(t_2-t_1)^{2H-\frac12},
\end{eqnarray}
where the last inequality follows as in (\ref{eq:n8.2}).

On the other hand, it is not difficult to see that
\begin{eqnarray}\label{eq:n22.1}
\lefteqn{\int_0^{t_1}\iou\lp\me_{t_2,s}(x,y)- \me_{t_1,s}(x,y)\rp^2dyds
}\nonumber\\
&&\le 2\int_{0}^{t_1}\iou\lp\int_{s+\ep}^{t_1+\ep}\lc
G_{t_2+2\ep-u}(x,y)-G_{t_1+2\ep-u}(x,y)\rc
\partial_u\kh(u,s)du\rp^2dyds\nonumber\\
&&+2 \int_{0}^{t_1}\iou\lp\int_{t_1+\ep}^{t_2+\ep}
G_{t_2+2\ep-u}(x,y)
\partial_u\kh(u,s)du\rp^2dyds\nonumber\\
&=&B_1+B_2.
\end{eqnarray}
Observe now that we can proceed as in (\ref{eq:n21.2}) to obtain 
\begin{equation}\label{eq:n22.2}
B_2\le c(t_2-t_1)^{2H-\frac12},
\end{equation}
and it is also readily checked that
\begin{eqnarray*}
B_1&\le& 2\int_{0}^{t_1+2\ep}\iou\lp\int_{s}^{t_1+2\ep}
|G_{t_2+2\ep-u}(x,y)-G_{t_1+2\ep-u}(x,y)|
\partial_u\kh(u,s)du\rp^2dyds\\
&=&
2H(2H-1)\int_{0}^{t_1+2\ep}\int_{0}^{t_1+2\ep}\iou {|u-v|^{2H-2}}
|G_{t_2+2\ep-u}(x,y)-G_{t_1+2\ep-u}(x,y)|\\
&&\qquad\quad\quad\quad\qquad\quad\cdot|G_{t_2+2\ep-v}(x,y)
-G_{t_1+2\ep-v}(x,y)|dydudv.
\end{eqnarray*}
Finally, the proof follows combining (\ref{eq:n8.3}), and 
(\ref{eq:n21.1})-(\ref{eq:n22.2}).

\end{proof}

\vspace{0.3cm}

Let us state now the main result of this paper.
\begin{theorem}\label{teo:ito}
Let $X$ be the process defined by  (\ref{sol:heateq}) and
 $f\in C_b^2(\R)$. Then, for
$t\in[0,T]$ and $x\in\ou$, the process $\ms_{t,x}f'(X)$
belongs to $\dom(\delta)$ and
$$
f(X(t,x))=f(0)+\mathbf{A_{1}}(t,x)+\mathbf{A_{2}}(t,x),
$$
where
$$
\mathbf{A_{1}}(t,x)=
\iot\iou \lp\ms_{t,x}f'(X)\rp(s,y) W(ds,dy)$$
and
$$
\mathbf{A_{2}}(t,x)=
\frac12 \iot f''(X(s,x))K_{x}(ds),
$$
with
$$
K_{x}(s)
=
H(2H-1) \int_0^s \int_0^s G_{2s-v_1-v_2}(x,x) |v_1-v_2|^{2H-2}
dv_1 dv_2.
$$
\end{theorem}

As mentioned before, in order to prove this theorem, we use the regularized
It\^o formula of Proposition \ref{pro:n2.6}, and we only need to
study the convergence of the terms $\mathbf{A_{1,\ep}}$ 
and $\mathbf{A_{2,\ep}}$ appearing there. However,
this analysis implies long and tedious calculations. 
This is why we have chosen to split the proof of our theorem
into a series of lemmas which will be given in the next section.

\subsection{Proof of the main result}
The purpose of this section is to present some technical results whose
combination provides us the proof of our It\^o's formula given at 
Theorem \ref{teo:ito}.
We begin with the convergence $\mathbf{A_{2,\ep}}\rightarrow
\mathbf{A_2}$, for which we provide first a series of lemmas.
\begin{lemma}\label{lem:n2.9}
Let $L_1^{\ep}$ be the function  defined on $\ott$ by
$$
L_1^{\ep}(s)=
\int_0^s dv_2\int_0^{v_2}dv_1\, 
G_{2(s+\ep)-v_1-v_2}(x,x)\kh(v_1+\ep,v_1)\partial_{v_2}\kh(v_2+\ep,v_1).
$$
Then $s\mapsto\partial_s L_1^\ep(s)$
converges to 0 in $L^1([0,T])$, as $\ep\downarrow 0$.
\end{lemma}
\begin{proof}
Note that by (\ref{dcp:ker}) we only need to study the convergence of
$\partial_s L_{11}^\ep(s)$, where
\begin{equation}\label{def:l11}
L_{11}^\ep(s)=\int_0^s dv_2\int_0^{v_2}dv_1\,
\frac{\kh(v_1+\ep,v_1)}{\sqrt{2(s+\ep)-v_1-v_2}}
\partial_{v_2}\kh(v_2+\ep,v_1) .
\end{equation}
Indeed, this term will show us the technique and the difficulties for
the remaining terms.

We will now proceed to a series of change of variables in
order to get rid of the parameter $s$ in the boundaries of
the integrals defining $L_{11}^\ep$:
using first the change of variable $z=\frac{v_2-v_1}{s-v_1}$, and
then $\theta=v_1/s$, we can write
\begin{multline*}
L_{11}^\ep(s)=c_Hs^{3-2H}\iou\lp\int_{s\theta}^{s\theta+\ep}
(u-s\theta)^{H-\frac32}u^{H-\frac12}du\rp\frac{(1-\theta)}{\theta^{2H-1}}\\
\cdot\iou\frac{(\ep +s\theta+zs(1-\theta))^{H-\frac12}}
{\sqrt{2\ep+s(1-\theta)(2-z)}}(zs(1-\theta)+\ep)^{H-\frac32}dzd\theta.
\end{multline*}
Hence, the change of variable $v=u-s\theta$ leads to
\begin{multline*}
L_{11}^\ep(s)=c_Hs^{3-2H}\iou\lp\int_{0}^{\ep}
v^{H-\frac32}(v+s\theta)^{H-\frac12}dv\rp\frac{(1-\theta)}{\theta^{2H-1}}\\
\cdot\iou\frac{(\ep +s\theta+zs(1-\theta))^{H-\frac12}}
{\sqrt{2\ep+s(1-\theta)(2-z)}}(zs(1-\theta)+\ep)^{H-\frac32}dzd\theta.
\end{multline*}
Therefore, by differentiating this expression in $s$, we end up with
a sum of the type
$$
\partial_sL_{11}^\ep(s)
=
\sum_{j=1}^{5}L_{11j}^\ep(s),
$$
where
\begin{multline*}
L_{111}^\ep(s)=
c_Hs^{2-2H}\iou\lp\int_{0}^{\ep}
v^{H-\frac32}(v+s\theta)^{H-\frac12}dv\rp\frac{(1-\theta)}{\theta^{2H-1}}\\
\cdot\iou\frac{(\ep +s\theta+zs(1-\theta))^{H-\frac12}}
{\sqrt{2\ep+s(1-\theta)(2-z)}}(zs(1-\theta)+\ep)^{H-\frac32}dzd\theta,
\end{multline*}
and where the terms $L_{112}^\ep,\ldots,L_{115}^\ep$, whose exact calculation
is left to the reader for sake of conciseness, are similar to
$L_{111}^\ep$.

Finally, we have
$$
L_{111}^\ep(s)\le c_H s^{-H}\lp\int_{0}^{\ep}v^{H-\frac32}dv\rp
\lp\iou\frac{(1-\theta)^{H-1}}{\theta^{2H-1}}d\theta\rp\iou z^{H-\frac32}dz
\le c_H \ep^{H-1/2} s^{-H},
$$
and it is easily checked that this last term converges to 0 in
$L^1(\ott)$. Furthermore, it can also be proved that 
$|L_{11j}^\ep(s)|\le cL_{111}^\ep(s)$ for $2\le j\le 5$, which
ends the proof.

\end{proof}

\begin{lemma}\label{lem:n2.10}
Let $L_{2}^\ep$ be the function defined on $\ott$ by 
$$
L_{2}^\ep(s)=
\int_0^s dv_2\int_0^{v_2}dv_1\, 
G_{2(s+\ep)-v_1-v_2}(x,x)
\partial_{v_1}\lp\int_{v_1}^{v_1+\ep}\kh(v_1+\ep,u)\partial_{v_2}
\kh(v_2+\ep,u)du\rp.
$$
Then $s\mapsto\partial_s L_2^\ep(s)$
converges to 0 in $L^1([0,T])$, as $\ep\downarrow 0$.
\end{lemma}

\begin{proof} As in the proof of Lemma \ref{lem:n2.9} we only show the 
convergence of $\partial_sL_{21}^\ep(s)$, where
$$
L_{21}^\ep(s)=
\int_0^s dv_2\int_0^{v_2}dv_1\, 
\frac{1}{\sqrt{2(s+\ep)-v_1-v_2}}
\partial_{v_1}\hat L(v_1,v_2),
$$
with
$$
\hat L(v_1,v_2)=
\int_{v_1}^{v_1+\ep}\kh(v_1+\ep,u)\partial_{v_2}
\kh(v_2+\ep,u)du.
$$
Towards this end, we will proceed again  to a series of changes
of variables in order to eliminate the parameter $s$ from the boundaries of
the integrals:
notice first that the definition of $\kh$, and
the change of variables $\theta=\frac{u-v_1}{r-v_1}$ and $z=r-v_1$
yield
\begin{multline*}
\hat L(v_1,v_2)=
c_H(v_2+\ep)^{H-\frac12}\int_0^\ep (v_1+z)^{H-\frac12}\iou
(v_1+\theta z)^{1-2H}\frac{z^{H-\frac12}}{(1-\theta)^{\frac32-H}}
\\
\cdot (v_2+\ep-v_1-\theta z)^{H-\frac32}d\theta dz.
\end{multline*}
Thus
\begin{align}\label{eq:n26.1}
&\partial_{v_1}\hat L(v_1,v_2)=
(v_2+\ep)^{H-\frac12}\Bigg[ c_H
\int_0^\ep (v_1+z)^{H-\frac32}\iou
(v_1+\theta z)^{1-2H}\frac{z^{H-\frac12}}{(1-\theta)^{\frac32-H}}\nonumber\\
&\hspace{10cm}\cdot (v_2+\ep-v_1-\theta z)^{H-\frac32}d\theta dz
\nonumber\\
&-c_H
\int_0^\ep (v_1+z)^{H-\frac12}\iou
(v_1+\theta z)^{-2H}\frac{z^{H-\frac12}}{(1-\theta)^{\frac32-H}}
(v_2+\ep-v_1-\theta z)^{H-\frac32}d\theta dz
\nonumber\\
& +c_H
\int_0^\ep (v_1+z)^{H-\frac12}\iou
(v_1+\theta z)^{1-2H}\frac{z^{H-\frac12}}{(1-\theta)^{\frac32-H}}
(v_2+\ep-v_1-\theta z)^{H-\frac52}d\theta dz
\Bigg].
\end{align}
Hence, it is easily seen that $L_{21}^\ep$ is a sum of terms of the form
\begin{eqnarray*}
Q_{\al,\beta,\nu}(s)&=&
\int_0^s dv_2\int_0^{v_2}dv_1\, 
\frac{(v_2+\ep)^{H-\frac12}}{\sqrt{2(s+\ep)-v_1-v_2}}
\int_0^\ep (v_1+z)^{\alpha}\\
&&\qquad\qquad\qquad\cdot 
\iou
(v_1+\theta z)^{\beta}\frac{z^{H-\frac12}}{(1-\theta)^{\frac32-H}}
(v_2+\ep-v_1-\theta z)^{\nu}d\theta dz\\
&=&s^2\iou d\eta\iou du\frac{\eta(s\eta+\ep)^{H-\frac12}}
{\sqrt{2(s+\ep)-us\eta-s\eta}}\int_0^\ep
(us\eta+z)^{\alpha}\\
&&\qquad\qquad\qquad\cdot 
\iou
(us\eta+\theta z)^{\beta}\frac{z^{H-\frac12}}{(1-\theta)^{\frac32-H}}
(s\eta+\ep-s\eta u+\theta z)^{\nu}d\theta dz,
\end{eqnarray*}
by applying the changes of variable $u=v_1/v_2$,
 and $\eta=v_2/s$. Differentiating this last relation, we are now 
able to compute $\partial_sL_{21}^\ep(s)$,
and see that this function goes to 0 as $\ep\downarrow 0$ 
in $L^1([0,T])$, similarly to what we did in the proof of 
Lemma \ref{lem:n2.9}.

\end{proof} 

\begin{lemma}\label{lem:n2.11}
Let $L_{3}^\ep$ be the function defined on $\ott$ by 
$$L_{3}^\ep(s)=H(2H-1)\int_0^s dv_2\int_0^{v_2}dv_1\, 
G_{2(s+\ep)-v_1-v_2}(x,x)(v_2-v_1)^{2H-2}.
$$
Then $\partial_sL_{3}^\ep(s)$ tends to $\frac{1}{2} K_x(ds)$ 
in $L^1([0,T])$, as $\ep\downarrow 0$.
\end{lemma}

\begin{proof} As in the proofs of Lemmas \ref{lem:n2.9} and \ref{lem:n2.10},
we only need to use the change of variables $z=v_1/v_2$ and 
$\theta=v_2/s$.

\end{proof}

\begin{lemma}\label{lem:n2.12} Let $X$ and $\xe$ be given in (\ref{rep:x})
and (\ref{def:xep}), respectively. Then $\xe(\cdot,x)$ converges to 
$X(\cdot,x)$ in $L^2(\Omega\times[0,T])$ and, for $t\in[0,T]$,
$\xe(t,x)$ goes to $X(t,x)$ in $L^2(\Omega)$, as $\ep\downarrow 0$.
\end{lemma}
\begin{proof}
 The result is an immediate consequence of the definitions of the 
processes $\xe(\cdot,x)$ and $X(\cdot,x)$, the fact that
$|\me_{t,s}(x,y)|\le c(t-s)^{H-1}s^{\frac12 -H}$ and of the dominated 
convergence theorem.

\end{proof}

We are now ready to study the convergence of the term
$\mathbf{A_{2,\ep}}$:
\begin{lemma}\label{lem:n2.13} Let $t\in[0,T]$ and
$x\in [0,1]$. Then the random variable
\begin{align*}
&\lefteqn{
B^\ep_2(t,x)}\\
&:=H(2H-1)\iot f''(\xe(s,x))\partial_s\lp\int_0^s
dv_2\int_0^{v_2}dv_1 G_{2(s+\ep)-v_1-v_2}(x,x)(v_2-v_1)^{2H-2}\rp ds
\end{align*}
converges to $\mathbf{A_2}(t,x)$ in $L^2(\Omega)$ as $\ep\downarrow 0$.
\end{lemma}
\begin{proof}
Since $f''$ is a bounded function, then
\begin{align*}
&\lefteqn{E\lp|B_2^\ep (t,x)-\mathbf{A_2}(t,x)|^2\rp}\\
&\le c\iot E\lp (f''(X(s,x))-f''(\xe(s,x)))^2\rp|\partial_s K_x(s)|ds\\
&+c \lp\iot|\partial_s K_x(s)-H(H-\frac12 )\partial_s
\int_0^sdv_2\int_0^{v_2}dv_1G_{2(s+\ep)-v_1-v_2}(x,x)(v_2-v_1)^{2H-2}|ds\rp^2
 .
\end{align*}
Hence
 the result is a consequence of Lemmas \ref{lem:n2.11} and \ref{lem:n2.12},
and the dominated convergence theorem.
\end{proof}

Now we study the convergence of $\mathbf{A_{1,\ep}}$ to $\mathbf{A_1}$ in
$L^2(\Omega)$.
\begin{lemma}\label{lem:n2.14}
 Let $X$ and $\xe$ be given in (\ref{rep:x})
and (\ref{def:xep}), respectively. Then, for $t\in[0,T]$ and $x\in\ou$,
$$E\lp\iot\iou\lc (\ms_{t,x}f'(X))(s,y)-(\mes_{t,x}f'(\xe))(s,y)\rc^2dyds\rp
\rightarrow 0$$
as $\ep\downarrow 0$.
\end{lemma}
\begin{proof}
We first note that
$$
E\lp\iot\iou\lc (\ms_{t,x}f'(X))(s,y)-(\mes_{t,x}f'(\xe))(s,y)\rc^2dyds\rp
$$
can be bounded from above by:
\begin{align}\label{eq:n29.1}
&cE\lp\iot\iou\lc\int_s^t(f'(X(r,x))-f'(X(s,x))\right.
\right.\nonumber\\
&\quad\quad\quad\quad\left.\left.\phantom{\iou}
-f'(\xe(r,x))+f'(\xe(s,x)))\partial_r
M_{r,s}(x,y)dr\rc^2dyds\rp\nonumber\\
&+c
E\lp\iot\iou\lc\int_s^t(f'(\xe(r,x))-f'(\xe(s,x)))(\partial_r
M_{r,s}(x,y)-\partial_r M_{r,s}^\ep (x,y))dr\rc^2dyds\rp\nonumber\\
&+ cE\lp\iot\iou\lc (f'(X(s,x))-f'(\xe(s,x)))M_{t,s}(x,y)\rc^2dyds\rp
\nonumber\\
&+ cE\lp\iot\iou\lc f'(\xe(s,x))(M_{t,s}(x,y)-M_{t,s}^\ep(x,y))\rc^2dyds\rp
\nonumber\\
&=c(B_1+\ldots+B_4).
\end{align}

Next observe  that
\begin{eqnarray*}
B_{2}&\le&\iot\iou\lc\int_s^t E\lp(
 f'(\xe(r,x))-f'(\xe(s,x)))^2\rp|\partial_r
M_{r,s}(x,y)-\partial_r M_{r,s}^\ep (x,y)|dr\rc\\
&&\quad\qquad\cdot\lc\int_s^t|\partial_\theta
M_{\theta,s}(x,y)-\partial_{\theta} M_{\theta,s}^\ep (x,y)|d\theta\rc dyds.
\end{eqnarray*}
Now notice that Proposition \ref{est:xep}  and
the inequality
\begin{eqnarray*}
\lefteqn{|\partial_r M_{r,s}^\ep (x,y)|}\\
&\le&c\lp\frac{r+\ep}{s}\rp^{H-\frac12}(r-s+\ep)^{H-2}
\lp\exp\lp-c_1\frac{(x-y)^2}{\ep+(r-s)}\rp+\exp\lp-c_1\frac{(x+y-2)^2}
{\ep+(r-s)}\rp\rp\end{eqnarray*}
imply, for $\beta$ small enough, that
$$E\lp(f'(\xe(r,x))-f'(\xe(s,x)))^2\rp|\partial_r
M_{r,s}(x,y)-\partial_r M_{r,s}^\ep (x,y)|$$
goes to 0 as $\ep\downarrow 0$ and that it is bounded by $cs^{\frac12-H}
(r-s)^{3H-\frac52-\beta}$. Thus 
\begin{equation}\label{eq:bonita}
B_{2}\rightarrow 0\end{equation}
 because of the
dominated convergence theorem.

Since $f'$ is a bounded function, then
$$B_4\le c\int_0^t\iou\lp M_{t,s}(x,y)-M_{t,s}^\ep(x,y)\rp^2dyds,$$
which goes to 0 due to the definition of $M^\ep$ and the dominated 
convergence theorem. Hence, by (\ref{eq:n29.1}),
 and (\ref{eq:bonita}), we only
need to show that $B_1+B_3\rightarrow 0$ as 
$\ep\downarrow 0$ to finish the proof.
This can been seen using Lemma \ref{lem:n2.12} 
and proceeding as the beginning of this proof.
\end{proof}
\begin{lemma}
Let $X$ and $\xe$ be given by (\ref{rep:x})
and (\ref{def:xep}), respectively. Then, for $t\in[0,T]$ and $x\in\ou$,
$\ms_{t,x}f'(X)$ belongs to Dom $(\delta)$. Moreover
$$\delta\lp \mes_{t,x}f'(\xe)\rp\rightarrow
\delta\lp\ms_{t,x}f'(X)\rp$$
as $\ep\downarrow 0$ in $L^2(\Omega)$.
\end{lemma}
\begin{proof}
 The result follows from Lemmas \ref{lem:n2.9}-\ref{lem:n2.14}, and from the
 fact that $\delta$ is a closed operator.
\end{proof}

\end{document}